\documentclass[11pt]{amsart}
\usepackage[latin1]{inputenc}
\usepackage{a4}
\usepackage{amssymb}
\usepackage{amsmath}
\begin{document}
\setlength{\parindent}{1.2em}
\def\noproof{{\unskip\nobreak\hfill\penalty50\hskip2em\hbox{}\nobreak\hfill%
       $\square$\parfillskip=0pt\finalhyphendemerits=0\par}\goodbreak}
\def\endproof{\noproof\bigskip}
\def\proof{\removelastskip\penalty55\medskip\noindent{\bf Proof. }}
\newtheorem{firstthm}{Proposition}
\newtheorem{thm}[firstthm]{Theorem}
\newtheorem{theorem}[firstthm]{Theorem}
\newtheorem{prop}[firstthm]{Proposition}
\vspace*{-2cm}

\title{A simple solution to Ulam's liar game with one lie}
\author{Deryk Osthus \and Rachel Watkinson}
\begin{abstract} \noindent
Ulam asked for the maximum number of questions required to determine an integer 
between 1 and $10^6$ by asking questions whose answer is `Yes' or `No' and where 
one untruthful answer is allowed.
Pelc showed that the number of questions required is 25.
Here we give a simple proof of this result.
\end{abstract}
\date{}
\maketitle \vspace{-.8cm}

\subsection*{Introduction}\label{intro}

We consider the following game between a questioner and a responder, first
proposed by Ulam~\cite{ulam}. 
(A variation of this game was independently proposed by 
R\'enyi, see~\cite{pelc1}.) 
The responder thinks of an integer $x \in \{1,\dots,n\}$ and the 
questioner must determine $x$ by asking questions whose answer is `Yes' or `No'.
The responder is allowed to lie at most $k$ times during the game.
Let $q_k(n)$ be the maximum number of questions needed by the questioner, under an optimal 
strategy, to determine $x$ under these rules.
In particular, Ulam asked for the value of $q_1(10^6)$ (as this is related to the 
well-known `twenty questions' game).
It follows from an observation of Berlekamp~\cite{berlekamp} that
$q_1(10^6) \ge 25$ and Rivest et al.~\cite{rivest} as well as 
Spencer~\cite{spencer1} gave bounds which imply that
$q_1(10^6) \le 26$.
Pelc~\cite{pelc2} was then able to determine $q_1(n)$ exactly for all $n$:
\begin{thm} \label{thmpelc}
\cite{pelc1} 
For even $n \in \mathbb{N}$, $q_1(n)$ is the smallest integer $q$ which satisfies $n \le 2^q/(q+1)$.
For odd $n \in \mathbb{N}$, $q_1(n)$ is the smallest integer $q$ which satisfies $n \le (2^q-q+1)/(q+1)$.
\end{thm}
In particular, his result shows that the lower bound of Berlekamp for $n=10^6$ was correct.
Shortly afterwards, 
Spencer~\cite{spencer1} determined $q_k(n)$ asymptotically (i.e.~for fixed $k$ and large $n$).
The values of $q_k(10^6)$ have been determined for all $k$.
These and many other related results are surveyed by Hill~\cite{hill},   
Pelc~\cite{pelc1} and Cicalese~\cite{cicalese}.
Here, we give a simple strategy and analysis for the game with at most one lie
which implies the above result of Pelc for many values of $n$.
\begin{thm}
If $n \le 2^\ell \le 2^q/(q+1)$ for some integer $\ell$, then
the questioner has a strategy which
identifies $x$ in $q$ questions if at most one lie is allowed. In particular, $q_1(n) \le q$.
\label{genstrat}\end{thm}
Below, we will give a self contained argument 
(Proposition~\ref{lowerprop}) which shows that
if $n$ also satisfies $n>2^{q-1}/q$, then the strategy in Theorem~\ref{genstrat} is optimal.
This implies that the bound in Theorem~\ref{genstrat} is optimal if $n=2^\ell$ for some 
$\ell \in \mathbb{N}$. 
More generally, Theorem~\ref{thmpelc} implies that for even $n$, Theorem~\ref{genstrat}
gives the correct bound if and only if we can find a binary power $2^\ell$ with $n \le 2^\ell \le 2^q/(q+1)$,
where $q$ is the smallest integer with $n \le 2^q/(q+1)$.
(Similarly, one can read off a more complicated condition for odd $n$ as well.)
In particular, if $n=10^6$, we obtain $q_1(10^6)=25$. To check this, note
that for $q=25$ and $\ell=20$, we have
$$
\lceil 2^{q-1} /q \rceil =671088 < n \le 1048576=2^{\ell} < 1290555= 
\lfloor 2^q/(q+1) \rfloor.
$$
If we compare the bounds from Theorems~\ref{thmpelc} and~\ref{genstrat}, 
then one can check that the smallest value where the latter gives a worse bound is $n=17$, 
where Theorem~\ref{genstrat} requires $9$ questions whereas $q_1(17)=8$. The smaller values are 
$q_1(2) = 3$, $q_1(3) =q_1(4)= 5$, $q_1(5)=\dots= q_1(8)=6$ and $q_1(9)=\dots=q_1(16)=7$.

More generally, it is easy to see that for any $n$ 
the strategy in Theorem~\ref{genstrat} uses at most two more 
question than an optimal strategy. Indeed, given $n$, let $\ell$ and $q$ be the smallest integers
satisfying $n \le 2^\ell \le 2^q/(q+1)$. So Theorem~\ref{genstrat} implies that $q$ questions suffice.
 Proposition~\ref{lowerprop} implies that
if $n>2^{q-3}/(q-2)$, then any successful strategy needs at least $q-2$ questions in the worst case.
To see that $n>2^{q-3}/(q-2)$, suppose that this is not the case. Then by assumption on $\ell$
we have $2^{\ell-1} <n \le 2^{q-3}/(q-2)$. So if $q \ge 4$ (which we may assume in view of the above discussion of
small values),
we have $2^\ell < 2^{q-2}/(q-2) \le 2^{q-1}/q$. This contradicts the choice of $q$.

Our proof of Theorem~\ref{genstrat} uses ideas from  Cicalese~\cite{cicalese} and Spencer~\cite{spencer}.
It gives a flavour of some techniques which are typical for the area.
Elsholtz (personal communication) has obtained another short proof for the case $n=10^6$.
Throughout, all logarithms are binary.

\medskip

From now on, we consider only the game in which at most one lie is allowed.
For the purposes of the analysis, it is convenient to allow the responder to 
play an adversarial strategy, i.e.~the responder does not have to think of 
the integer $x$ in advance (but does answer the questions so that there always is at least
one integer $x$ which fits all but at most one of the previous answers). 
The questioner has then determined $x$ as soon as
there is exactly one integer which fits all but at most one of the previous answers.
We analyze the game by associating a sequence of states $(a,b)$ with the game.
The state is updated after each answer.
$a$ is always the number of integers which fit all previous answers and $b$ is the number of integers which 
fit all but exactly one answer.
So initially, $a=n$ and $b=0$.
The questioner has won as soon as $a+b \le 1$.
If there are $j$ questions remaining in the game and the state is $(a,b)$, then
we associate a weight $w_j(a,b):=(j+1)a+b$ with this state.
Also, we call the integers which fit all but one exactly answer \emph{pennies}
(note that each of these contributes exactly one to the weight of the state).

For completeness, we now give a proof of the lower bound mentioned in the introduction.
As mentioned above, the fact is due to Berlekamp~\cite{berlekamp}, 
see also~\cite{cicalese,pelc2,rivest} for the argument.
The proof has a very elegant
probabilistic formulation which generalizes more easily to the case of $k \ge 1$ lies
(see Spencer~\cite{spencer}).
\begin{prop} \label{lowerprop}
If $n>2^{q-1}/q$, 
then the questioner
does not have a strategy which determines $x$ with $q-1$ questions.
\end{prop}
\proof
Note that our assumption implies that the initial weight satisfies $w_{q-1}(n,0)>2^{q-1}$.
It is easy to check that before each answer, the sum of the weights of the two possible new states
$(a_{yes},b_{yes})$ and $(a_{no},b_{no})$
is equal to the weight of the current state $(a,b)$, i.e.
\begin{equation} \label{sum}
w_j(a,b)=w_{j-1}(a_{yes},b_{yes})+w_{j-1}(a_{no},b_{no}).
\end{equation}
To see this, observe that $a=a_{yes}+a_{no}$ and $a+b=b_{yes}+b_{no}$ and substitute this into the
definition of the weight functions.
(\ref{sum}) implies that
the responder can always ensure that the new state $(a',b')$ (with $j$ questions remaining)
satisfies 
\begin{equation} \label{half}
w_j(a',b') \ge w_{j+1}(a,b)/2 \ge w_{q-1}(n,0)2^{-(q-1-j)}>2^j.
\end{equation}
Thus responder can ensure that the final state has weight greater than one.
We also claim that this game never goes into state $(1,0)$.
(Together, this implies that the final state consists of more than one penny, which means that
the responder wins).
To prove the claim, suppose that we are in state $(1,0)$ with $j-1$ questions to remaining. 
Then the previous state must have been $(1,t)$ for some $t>0$.
Note that~(\ref{half}) implies that $w_j(1,t)> 2^j$.
On the other hand, the assumption on the strategy of the responder implies that
$w_{j-1}(1,0) \ge w_{j-1}(0,t)$.
Combined with~(\ref{sum}), this means that 
$w_j(1,t)= w_{j-1}(1,0) + w_{j-1}(0,t) \le 2w_{j-1}(1,0)=2j$.
But $2j <2^j$ has no solution for $j \ge 1$, and so we have a contradiction.
\endproof

\subsection*{Proof of Theorem~\ref{genstrat}}
Note that the weight of the
initial state is
$w_q(n,0) = n(q+1) \leq 2^q.$ 
By making $n$ larger if necessary, note that we may assume that $\log n= \ell$, 
for some $\ell \in \mathbb{N}$.
So $\ell \leq q - \log (q+1)$.
Since $\ell \in \mathbb{N}$, this implies
\begin{align}{\ell \leq q - \lceil\log (q+1)\rceil \label{lv}}.\end{align}
 Consider each integer 
\begin{math}{n = 2^\ell}\end{math} in its binary form, i.e.~we have $2^\ell$
strings of length $\ell$. The questioner performs a binary search on
these numbers by asking questions of the form `Is the value of $x$
in position $i$ a 1?'. The binary search on the search space $\{1,\dots,n\}$
uses exactly $\ell$ questions and as a result we obtain $\ell+1$ possible
binary numbers for $x$. There is exactly one integer which satisfies
all the answers. There
are also $\ell$ integers which satisfy all but one answer.
Therefore, after the binary search has been performed we are in state $(1,\ell)$.
Moreover, $w_{q-\ell}(1,\ell)  = 1\cdot (q-\ell+1) + \ell\cdot 1 =
q+1.$

Let $p=q-\ell$. By (\ref{lv}), it now suffices to identify $x$ within
\begin{math}{p :=\lceil\log (q+1)\rceil }\end{math} questions. Note
that the weight of the state satisfies \begin{math}{2^{p-1}<
\;w_{q-\ell}(1,\ell)\; \leq 2^p}\end{math}. 
Suppose that $q+1$ is not
a power of 2. It
is easy to see that we can add pennies to the state until the total
weight is equal to $2^p$, as the addition of pennies will only make
the game harder for the questioner. Suppose that we now have $r$
pennies in total, so we obtain the new state $P^* =
(1,r)$, with $r \geq \ell$, where the weight of $P^*$ equals
$2^p$.
Thus
 \begin{align}{p + 1+r = w_p(1,r) =
2^p.}\label{**}\end{align}
 We now have two cases to consider:

\vspace{0.2cm}

\noindent\textbf{Case One:} If $ r< p+1$, then~(\ref{**}) implies that 
$p+1  > 2^{p-1}$, which holds if and
only if $p\leq 2$. This means that we have one nonpenny and at most
two pennies. It is easy to see that the Questioner can easily identify $x$ using
two more questions in this case.

\vspace{0.2cm}

\noindent\textbf{Case Two:} Suppose $r\geq p+1$. This implies that $2^{p-1} \geq p+1$ and thus $p>2$.
We know that the total weight of this state is even and so we wish to
find a set, say $A_p$, such that when a question is asked about it,
regardless of the responder's reply, the weight is exactly halved.
Assume that $A_p$ contains the nonpenny and $y$ pennies and that the
weight of $A_p$ is equal to $2^{p-1}$. Suppose that the answer to `Is $x \in A_p$?' is `Yes'. 
Then the weight of the resulting state is $p + y$ 
(since we are left with one nonpenny of weight $p$ and $y$ pennies).
If the answer is `No', the resulting state has weight $r+1-y$
(since the nonpenny has turned into a penny and the $y$ pennies have been excluded).
Thus we wish to solve $r+1- y = p +y$, which gives
\begin{equation} \label{y}
 y = \frac{1}{2}(r+1-p).
\end{equation}
Note also that
that (\ref{**}) implies $r+1-p$ is even and
so $y$ is an integer. Moreover, the condition
$r \geq p+1$ implies that $ y \geq 1$.

So suppose that the questioner chooses $A_p$ as above and asks `Is $x\in
A_p$?'. 
If the responder replies `Yes', we obtain a position
$P'$, which consists of one nonpenny and $y$ pennies, i.e.~$P' =(1,y)$, 
which has weight $2^{p-1}$. If $p-1 = 2$, then by Case~1,
the questioner can easily identify $x$. If $p-1 > 2$, we redefine
$r$ such that $r :=y$ and then calculate the new value of $y$ by
(\ref{y}), to obtain a new set $A_{p-1}$. The questioner continues 
inductively with $A_{p-1}$ instead of $A_p$, so the next question will 
be `Is $x\in A_{p-1}$?'. 
If the responder replies `No' to the original question `Is $x\in
A_p$?' then we obtain a position
$P'$ which consists only of pennies, i.e.~$P'=(0,r-y+1)$. 
Again, this has weight $2^{p-1}$. Since we have
$p-1$ questions remaining we perform a binary search on the $r-y-1=2^{p-1}$
pennies remaining and after $p-1$ questions we will have identified
$x$.
 
Note that eventually, the answer to the question `Is $x\in A_i$?'
must be either `No' or it is `Yes' and we have $i-1=2$ as well as a new weight of $2^{i-1}$
(in which case
there are 2 questions and at most one nonpenny and two pennies remaining). 
By the above arguments, the questioner
can find the integer $x$ in the required total number $q$ of questions in both cases,
which completes the proof of the theorem.\endproof

In case $n=10^6$, the above strategy would mean that after $20$ questions, we would be in 
state $(1,20)$ and have weight $w_5(1,20)=26$. Our aim is to find $x$ within $5$ more 
questions. We add $6$ pennies to obtain the state $(1,r)$ with $r=26$ and weight $2^p$,
where $p=5$. Thus~(\ref{y}) gives $y=11$. So $A_5$ consists of the nonpenny and $11$
pennies. If the answer is `Yes', then $A_4$ consists of the nonpenny and $4$ pennies.
If the answer is `No', we have 16 pennies left and can find $x$ after $4$ more questions
by using binary search.

\enlargethispage{1cm}

\medskip

\enlargethispage{0.6cm}
{\footnotesize \obeylines \parindent=0pt
Deryk Osthus \& Rachel Watkinson
School of Mathematics, University of Birmingham, Edgbaston, Birmingham, B15 2TT, UK
}

{\footnotesize \parindent=0pt
\it{E-mail address}:
\tt{osthus@maths.bham.ac.uk, rachel.watkinson@btinternet.com}}
\end{document}